\documentclass{amsproc}

\usepackage{a4wide}
\usepackage{amsmath}
\usepackage{enumerate}

\usepackage{amsmath,amsthm}
\usepackage{amsfonts}
\usepackage{amssymb}
\usepackage{graphicx, color, soul}
\usepackage[T1]{fontenc}
\usepackage{verbatim}
{\bf}{\rm}{\rm}

\newtheorem*{theorem*}{Theorem}
\newtheorem*{cor*}{Corollary}

\newtheorem{theorem}{Theorem}
\newtheorem{cor}{Corollary}
\newtheorem{prop}{Proposition}

\newtheorem{lem}{Lemma}

{\rm}{\rm}

\def\Real{\mathbb{R}}

\def\ad{\text{\rm ad}}

\def\SO{\text{\rm SO}}

\def\CO{\text{\rm CO}}

\newcommand{\frakk}{\mathfrak{k}}
\newcommand{\g}{\mathfrak{g}}
\newcommand{\f}{\mathfrak{f}}
\newcommand{\heis}{\mathfrak{heis}}

\newcommand{\su}{\mathfrak{su}}

\newcommand{\h}{\mathfrak{h}}

\newcommand{\isom}{\mathfrak{isom}}

\newcommand{\conf}{\mathfrak{conf}}

\renewcommand{\so}{\mathfrak{so}}
\newcommand{\co}{\mathfrak{co}}
\newcommand{\hol}{\mathfrak{hol}}

\newcommand{\zl}{\rtimes}

\def\id{\mathop\text{\rm id}\nolimits}

\def\ad{\mathop\text{\rm ad}\nolimits}

\def\tr{\mathop\text{\rm tr}\nolimits}
\def\Ric{\mathop\text{\rm Ric}\nolimits}

\def\spa{\mathop\text{{\rm span}}\nolimits}

\usepackage{stackrel}

\newcommand{\be}{\begin{equation}}
\newcommand{\ee}{\end{equation}}

\let\geq=\geqslant

\begin{document}
	
	\title{Conformally homogeneous Lorentzian spaces}

	\author{Dmitri V. Alekseevsky}\thanks{$^1$Higher School of Modern Mathematics MIPT, 		1 Klimentovskiy per., Moscow, Russia}
	
	\author{Anton S. Galaev}\thanks{$^2$University of Hradec Kr\'alov\'e, Faculty of Science, Department of Mathematics, Rokitansk\'eho 62, 500~03 Hradec Kr\'alov\'e,  Czech Republic,
		E-mail: anton.galaev(at)uhk.cz}

	\begin{abstract}
		
	We prove that if a 1-connected non-conformally flat conformal Lorentzian manifold $(M,c)$  admits a connected essential  transitive group of conformal transformations, then there exists a metric $g\in c$ such that $(M,g)$ is a complete  homogeneous plane wave. This completes the classification of 1-connected Lorentzian manifolds that admit a transitive essential conformal group.	We also prove that the
		group of conformal transformations of
		a non-conformally flat 1-connected homogeneous plane wave $(M,g)$ consists of homotheties,  and it is a 1-dimensional extension of the group of isometries.

		\vskip0.1cm
		
		{\bf Keywords}: conformal manifold; homogeneous space; conformal transformation; plane wave; Lorentzian manifold
		\vskip0.1cm
		
		{\bf AMS Mathematics Subject Classification 2020:}    53C18; 53B30; 53C50; 53C30
		
		
	\end{abstract}

	\maketitle

\section{Introduction}

It is well known that any Riemannian manifold which admits an essential group of conformal transformations is conformally equivalent to the standard
sphere or the Euclidean space. It is the Lichnerowicz conjecture, proved in  \cite{A,A2,Fer,Obata}.
On the other hand, there are many examples of pseudo-Riemannian (in particular
Lorentzian) manifolds with essential conformal group. Frances \cite{F,F1} constructed the first examples of conformally essential compact Lorentzian manifolds.
Podoksenov \cite{P} found examples of essential conformally homogeneous
Lorentzian manifolds. A local description of Lorentzian manifolds with
essential group of homotheties was given in \cite{A1}.

We study  essential conformally
homogeneous conformal Lorentzian manifolds $(M=G/H, c)$, i.e.,   conformal manifolds  with   transitive   group $G$  of  conformal transformations
which   does not   preserve any metric from  the   conformal  class~$c$.
As in \cite{A2017},  such conformal  manifolds $(M=G/H,c)$ may be split into two types:

\begin{itemize}
	\item[\bf A.]  Manifolds   with non-faithful   isotropy  representation
$$ j: \mathfrak{h} \to \mathfrak{co}(V),\quad V = \mathfrak{g}/\mathfrak{h}= T_oM $$
of  the   stability   subalgebra  $\mathfrak{h}$.

\item[\bf B.]  Manifolds   with   faithful  isotropy  representation $j$.
\end{itemize}

In \cite{A2017}, Lorentzian manifolds of type A were classified. We review this result in Section \ref{secA}.
In particular, manifolds of type A are conformally flat.
In this paper we classify 1-connected non-conformally flat essentially conformally homogeneous Lorentzian manifolds. These homogeneous spaces are of type B and  they  are exhausted by  the homogeneous  plane waves. More precisely
we prove the following Main  Theorem.

\begin{theorem}\label{ThMain}
Let $(M=G/H, c)$ be a   1-connected non-conformally flat  essential homogeneous  conformal Lorentzian manifold.  
Then there exists a metric $g\in c$ such that $(M,g)$ is a complete homogeneous plane wave, and   a  transitive  subgroup $G’ \subset   G$   preserves the metric $g$. 
\end{theorem}

We say      that   a  Lorentzian manifold $(M,g)$  is   homogeneous  (resp.,  conformally  homogeneous) if the  isometry  group  (resp.,  the  conformal  group) acts transitively on $M$. The local form of homogeneous plane waves was found in \cite{BO}. Recently in \cite{HMZ} it was shown that a 1-connected
homogeneous plane wave admits global Brinkmann coordinates. The Lie algebra of conformal Killing vector fields  of these spaces are known. We recall these results in Section \ref{secHpw}.

Section \ref{secpMTh}  contains the proof of the Main Theorem \ref{ThMain}. The proof consists of three steps. 
As the first step, we prove that the isotropy subalgebra $j(\h)\subset\co(V)$ contains an element $D$ of a particular form. In the second step  we show that the transitive conformal group  $G$  of the manifold $(M,c)$ contains a  Lie subgroup $F$ which has an open orbit $U\subset M$ and acts on $U$ by isometries of the restriction $g_U$ of a  metric $g \in c$  from  the conformal class. We prove that $(U,g_U)$ is an $F$-homogeneous plane wave. In the third step we consider the embedding of the Lie algebra of $G$  into the Lie algebra  of conformal Killing vector fields of the  homogeneous plane wave $(U,g_U)$. This implies that $F\subset G$ is a normal subgroup, and consequently $U=M$. This allows  to complete the proof of  Theorem \ref{ThMain}.

In the  recent work \cite{HSII2024} it is  proved that the Lie algebra of conformal Killing vector fields of a 1-connected homogeneous plane wave is a 1-dimensional extension of the Lie algebra of  Killing vector fields  and it consists of homothetic vector fields. 
We extend this result and we 
compute the conformal group of 1-connected  homogeneous plane waves.

\begin{theorem}\label{ThConf}
	Let $(M=G/H, g)$ be a 1-connected non-conformally flat Lorentzian essential   conformally homogeneous  manifold    of  dimension at least four. Then the conformal   group  $\mathrm{Conf}(M,g)$  consists of homotheties and it is a 1-dimensional extension of the group of isometries.
	\end{theorem}

The proof of Theorem \ref{ThConf} is given in Section \ref{secCGHpvGroup}. 	Let $(M=G/H, g)$ be a 1-connected non-conformally flat Lorentzian essential  conformally homogeneous  manifold    of  dimension at least four. According  to Theorem  \ref{ThMain},    any   such manifold   is  conformally  diffeomorphic    to  a plane  wave manifold.  Each conformal diffeomorphism $a$ of $(M,g)$ preserves the conformal Weyl curvature tensor. Analyzing the Ricci tensor and the Weyl tensor of $(M,g)$, we see that $a$ preserves also the Ricci tensor of $(M,g)$, i.e., $a$ is a Liouville transformation of $(M,g)$ in the sense of \cite{KRLeovTrans}.   A result from \cite{KRLeovTrans} shows that $a$ is a homothety transformation.   

In Section \ref{secEx} we consider two special cases. 
 First,
the case of  1-connected  Lie groups  $G$ with left-invariant Lorentzian metrics $g$ that admit essential conformal transformations   induced by derivations of the corresponding Lie algebras.    
We give a new proof of  the   classification of  such   Lorentzian manifolds  $(G,g)$,  obtained   in \cite{ZC2021,ZC2024}.  
Second, we give  a criterion for a homogeneous plane wave  to admit a simply  transitive Lie  group  of isometries.

\vskip0.3cm

{\bf Acknowledgements.}
The authors are thankful to Vicente Cortés, Thomas Leistner and  the anonymous referee for valuable comments and suggestions
that have significantly improved the manuscript, in particular, the proof of the main theorem. 
The research of D.A. was supported by Basis-foundation-Leader grant  no. 22-7-1-34-1 and by the MSHE "Priority 2030" strategic academic leadership program.  A.G. was supported by the project GF24-10031K of Czech Science Foundation (GA\v{C}R).

\section{Conformally homogeneous spaces with non-faithful isotropy representation}\label{secA}

In this section we shortly discuss  the results from \cite{A2017} that give a description of all homogeneous Lorentzian  manifolds of type A.
In particular we will see that all these manifolds are conformally flat.

Let $(M=G/H, c)$ be  a  conformally  homogeneous
pseudo-Riemannian manifold  of  signature  $$(k, \ell) =
(-\cdots-,+\cdots +).$$ Denote by $\g$ and $\h$ the Lie algebras of the Lie groups $G$ and $H$, respectively.
 Let $$j^H:H\to \text{\rm CO}(V),\quad  j:\mathfrak{h}
\to \mathfrak{co}(V)$$ be  the isotropy representations of the stability
subgroup $H$ and the stability subalgebra $\mathfrak{h}$ of the
point $o=eH \in M$ in the tangent space $V =T_oM$.
Since, by the assumption, the kernel $\ker j$ of the representation $j$ of $\h$ in $\co(V)$ is non-trivial.
There is a filtration 
$$\mathfrak{g}_{-1}= \mathfrak{g} \supset
\mathfrak{g}_{0}= \mathfrak{h}\supset \mathfrak{g}_{1}= \mathrm{ker}j \supset
\mathfrak{g}_{2}=0.$$
The  associated $|1|$-graded Lie algebra is
$$ \mathrm{gr}(\mathfrak{g})=
\mathfrak{g}^{-1} \oplus \mathfrak{g}^0 \oplus \mathfrak{g}^1  =  V \oplus \mathfrak{g}^0 \oplus
\mathfrak{g}^1, $$
where $$V = \mathfrak{g}/\mathfrak{h},\quad  \mathfrak{g}^0 =
\mathfrak{h}/\mathfrak{g}_{1} = j(\mathfrak{h})\subset\co(V), \quad
\mathfrak{g}^1 = \mathfrak{g}_1 = \mathrm{ker}j.$$ Here and in what follows, the direct sum symbol $\oplus$ refers to a direct sum of vector spaces.
Since $\mathrm{gr}(\g)$ is a $|1|$-graded Lie algebra, the space $\g^1$ is contained in the first prolongation $(\g^0)^{(1)}$ of $\g^0\subset\co(V)$. Recall that 
$$(\g^0)^{(1)}=\{\varphi\in V^*\otimes \g^0,\,\varphi(X)Y=\varphi(Y)X,\,\forall X,Y\in V \}.$$
Recall also that $(\co(V))^{(1)}\cong V^*$. Thus, $\g^1$ may be identified with a $\g^0$-invariant subspace of $V^*$.

Let us  describe the standard model of  conformally  flat   pseudo-Riemannian  conformal manifold. 
The  projectivisation  $$S^{k, \ell} =P \mathbb{R}^{k+1, \ell+1}_0 \subset P\mathbb{R}^{k+1, \ell+1} $$ of  the isotropic   cone  $\mathbb{R}^{k+1, \ell+1}_0 \subset \mathbb{R}^{k+1, \ell+1}  $
carries   a conformally   flat  conformal    structure  of   signature  $(k, \ell)$. 
The orthogonal group $ \mathrm{SO}(k + 1, \ell+ 1)$  acts  transitively on  
$ S^{k, \ell} $     and $ S^{k, \ell} $ is  the    maximally   homogeneous space  represented by 
$$\mathrm{SO}(k + 1, \ell+ 1)/H,$$
where  $H$ is    the stability  subgroup  isomorphic   to  the
group   of  similarities     $$\mathrm{Sim}(V)=  \text{\rm CO}(V)\cdot V$$ of  the pseudo-Euclidean vector  space  $ V= \mathbb{R}^{k,\ell}$.
The associated  graded  Lie   algebra is
$$
\mathrm{gr}(\mathfrak{so}(k+1,\ell+1))\cong
\mathfrak{so}(k+1,\ell+1)= V \oplus \mathfrak{co}(V) \oplus V^*,
$$
where $V^* =\mathfrak{co}(V)^{(1)}$
is  the  first prolongation   of   $\mathfrak{co}(V)$.

Let    $$\mathfrak{g} = \mathfrak{g}^{-1} \oplus \mathfrak{g}^0 \oplus \mathfrak{g}^1 = V \oplus \mathfrak{g}^0 \oplus \mathfrak{g}^1 $$
be  a $|1|$-graded Lie algebra with $\g^0\subset\co(V)$ and $\g^1\neq 0$. It is obvious that $\mathfrak{g}$ may be considered as a  subalgebra  of the   graded  Lie   algebra
$\mathfrak{so}(k+1, \ell+1)$. Denote     by  $\hat G$   the    corresponding connected  Lie  subgroup of $\SO(k+1,\ell+1)$   and    by  $\hat H$  the   connected  subgroup generated   by  the subalgebra $\mathfrak{h}= \mathfrak{g}^0 + \mathfrak{g}^1\subset\g$. Since $\g$ contains $V$,  the orbit
$$\hat Go=\hat G/\hat H\subset S^{k,\ell}$$
is open, and it admits the canonical flat conformal structure. 

Let    $(M= G/H,c)$  be a   conformally  homogeneous  manifold of type A.   The associated $|1|$-graded Lie  algebra  $\mathrm{gr}(\mathfrak{g})$ has a  natural
embedding  into  the  $|1|$-graded Lie  algebra  $\mathfrak{so}(k+1, \ell+1)$ as a graded subalgebra.
Suppose that $\g$ is isomorphic 
to  the associated    graded  Lie  algebra $\mathrm{gr}(\mathfrak{g})$. 
Let $\hat G$ and $\hat H$ be the Lie groups associated to $\mathrm{gr}(\g)\cong\g$  as above.  Then $(M= G/H,c)$ is locally conformally diffeomorphic to $\hat G/\hat H$ with the canonical flat conformal structure. This implies

\begin{theorem}\cite{A2017} Let    $(M= G/H,c)$ be a   conformally  homogeneous  manifold of type A. Suppose that $\g$ is isomorphic 
		 to  the associated    graded  Lie  algebra $\mathrm{gr}(\mathfrak{g})$.
		 Then the manifold $(M= G/H,c)$ is locally conformally diffeomorphic to $S^{k,\ell}$ with the canonical flat conformal structure. In particular, $(M= G/H,c)$ is   conformally  flat.
\end{theorem}

Let    $(M= G/H,c)$ be a   conformally  homogeneous  manifold of type A. Suppose that $\g^0=\co(V)$. Then $\g^1=V^*$. This implies that $\g\cong \mathfrak{so}(k+1, \ell+1)$, i.e., $\mathrm{gr}(\g)\cong \g$. Hence the manifold is isomorphic to the standard model. It remains to consider the case when $\g^0$ is a proper subalgebra of $\co(V)$ and $\g$ is not isomorphic 
to  $\mathrm{gr}(\mathfrak{g})$. In this case the  following  result holds

\begin{theorem}\cite{A2017} \label{ThF}
	 Let  $(M=G/H,c)$ be  a conformally  homogeneous    Lorentzian
	manifold of  type  A  such that the  isotropy  algebra
	$j(\mathfrak{h})$ is  a proper  subalgebra of $\mathfrak{co}(V)$. If  the Lie algebra  $\mathfrak{g}$ is  not isomorphic  to  the associated    graded  Lie  algebra $\mathrm{gr}(\mathfrak{g})$,  then $M$ is  conformally diffeomorphic  to   the Fefferman
	space.
\end{theorem}

Recall that the Fefferman space of the Lorentzian signature $(1,2m+1)$  is  defined  as   the  manifold $F$ of  real
isotropic lines in $\mathbb{C}^{1,m+1}$. The Fefferman space is the homogeneous manifold $\mathrm {SU}(1,m+1)/H,$
where $H\subset \mathrm{SU}(1,m+1)$ is the stabilizer of a real isotropic line. In \cite{A2017} it is shown that the value of the  curvature tensor  of the Fefferman space at each point coincides with the value of the curvature of a conformally flat Cahen-Wallach symmetric space. In particular, the Fefferman space is conformally flat.

The idea of the proof of Theorem \ref{ThF} is the following. 
 The starting point is a construction of  a  special  element $D$ of the first prolongation   $(\g^0)^{(1)}$ of  the Lie  algebra  $\g^0$ which defines a  |2|-grading  of the Lie  algebra $\g$.  Analyzing the  Jacobi identity, one may check   that the graded Lie  algebra  $\g$ is isomorphic to the   Lie algebra    $\su(1, m + 1)$   with the canonical |2|-grading.   This  implies the  Theorem.

\begin{cor}\label{Corconfflat}
	If $(M,c)$ as a conformally homogeneous spaces of type A, then $(M,c)$ is conformally flat.
\end{cor}

\section{Homogeneous plane waves}\label{secHpw}

Recall that a Lorentzian manifold $(M,g)$ is called a plane wave if
there exists a parallel isotropic vector field $p$ such  that  the curvature tensor $R$ of $(M,g)$ satisfies the conditions
\begin{equation}\label{nablaXR} R(X,Y)=0,\quad
\nabla_X R=0\quad \text{for all vector fields  $X,Y$ orthogonal to $p$}.
\end{equation}
The metric $g$  of  a plane wave may be written locally in the form
\begin{equation}\label{localplanewave}
g=2dvdu+\sum_{i=1}^n (dx^i)^2+a_{ij}(u)x^ix^j(du)^2,
\end{equation}
where  $a_{ij}(u)$ is  a symmetric matrix of functions. The metric \eqref{localplanewave} is conformally flat if and only if
$$a_{ij}(u)=\delta_{ij}b(u),$$ where $b(u)$ is a function.
Recently it was shown in \cite{HMZ}, using the results from \cite{BO},  that  a 1-connected homogeneous non-flat plane wave $(M,g)$ is globally isometric to  one of the following  model spaces:

\begin{itemize}
	\item[\bf (a)] the space $\Real^{n+2} = \Real \times \Real^n\times \Real$
with the metric
$$g=2dvdu+\sum_{i=1}^n (dx^i)^2+\left(e^{uF}Be^{-uF}\right)_{ij}x^ix^j(du)^2,$$
	\item[\bf (b)] 	 the space $\Real\times\Real^n\times\Real_{>0}$ with the metric
	$$g=2dvdu+\sum_{i=1}^n (dx^i)^2+\left(e^{\ln(u)F}Be^{-\ln(u)F}\right)_{ij}x^ix^j\frac{(du)^2}{u^2}.$$
\end{itemize}
Here $B$ and $F$ are respectively symmetric and skew-symmetric matrices. The metrics of type (a) are geodesically complete, while the metrics of type (b) are not geodesically complete.

Note that each homogeneous plane wave of type (b) is globally conformally diffeomorphic to  a homogeneous plane wave of type (a), see, \cite{HSII2024}. Indeed, the  coordinates transformation  $$v\mapsto v-\frac{1}{4}\sum_{i=1}^n(x^i)^2,\quad x^i\mapsto e^{\frac{u}{2}}x^i,\quad u\mapsto e^u,$$   transforms  the  metric (b)  into the metric 
$$g=e^u\left(2dvdu+\sum_{i=1}^n (dx^i)^2+\left(e^{uF}\left(B-\frac{1}{4}\id\right)e^{-uF}\right)_{ij}x^ix^j(du)^2\right).$$
\smallskip

In order to describe the structure of the Lie algebra of Killing vector fields of a homogeneous plane wave, we introduce some notation that we will use throughout this paper.
Denote by $V$  the Minkowski space $\Real^{1,n+1}$ with the metric $(\cdot,\cdot)$. We identify  the Lorentz Lie algebra $\so(V)$ with the space of bivectors $\wedge^2 V$ in such a way that
$$(X\wedge Y)Z=(X,Z)Y-(Y,Z)X,\quad\forall X,Y,Z\in V.$$
Let $p,e_1,\dots,e_n,q$  be a Witt basis of $V$. We denote by $E$ the Euclidean space spanned by the vectors $e_1,\dots,e_n$.
We obtain the decomposition 
\begin{equation}\label{decompsoVC}\so(V)=\Real p\wedge q\oplus\so(E)\oplus p\wedge E\oplus q\wedge E,\end{equation} where 
$\co(E)=\Real p\wedge q+\so(E)\subset\so(V)$ is the maximal reductive subalgebra, and $p\wedge E$ and $q\wedge E$ are $\ad_{\co(E)}$-invariant commutative subalgebras.

The Heisenberg Lie algebra  may be defined as the  Lie  algebra
$$     \heis(E)= p\wedge E\oplus E\oplus \Real p$$  with the only  non-zero  Lie bracket
$$[p\wedge Y,X]=(p\wedge Y)X=-(Y,X)p.$$
The  orthogonal Lie algebra $\so(E)$ acts on $\heis(E)$ in the obvious way. 

According to \cite{HMZ}, the following theorem holds true.

\begin{theorem}\label{ThIsotrplw}
Let $(M,g)$ be a  1-connected homogeneous plane wave of type (a) or (b) as above and  $\frakk\subset\so(E)$  the subalgebra commuting with $B$ and~$F$. Then the  isometry Lie  algebra  and  the  stability  subalgebra of $(M,g)$  are  given by   
\begin{equation}\label{algisom}\isom(M,g)=(\Real q\oplus\frakk)\zl\heis(E)=\frakk\oplus p\wedge E\oplus V,
\end{equation} 
$$\isom(M,g)_o=\frakk+p\wedge E,$$ 
where the Lie algebra $\frakk+p\wedge E$ acts on $V$ in the standard way and  
\begin{align*}
[q,p]&=\lambda p,\quad [p,X]=0,\quad [X,Y]=0,\\
[q,p\wedge X]&=p\wedge(\lambda\id_E +F)X-X,\\
[q,X]&=p\wedge BX+FX,
\end{align*}
for all $X,Y\in E$. Here $\lambda=0$ for the spaces of type (a), and $\lambda=1$ for the spaces of type (b).
\end{theorem}

Let us now consider the Lie algebra of conformal Killing vector fields of a homogeneous plane wave.
It is clear that, for any non-zero $\lambda\in\Real$, the transformation
\begin{equation}\label{homoth}
(v,x^i,u)\mapsto (\lambda^2v,\lambda x^i,u)
\end{equation} is a homothety transformation of arbitrary plane wave metric, see, e.g., \cite{Blau,HSII2024}. This 1-parameter Lie group defines the 
conformal Killing vector field $$D=2v\partial_v+x^i\partial_{x^i}.$$
From \cite[Corollary 2]{HSII2024} it follows that
$$[D,p]=2p,\quad [D,X]=X,\quad X\in E,\quad [D,q]=0.$$
This shows that $\ad_D$ acts on $V=\isom(M,g)/\isom(M,g)_o$ as the 
endomorphism $\id_V-p\wedge q.$  The following theorem is proved in \cite{HSII2024}. 

\begin{theorem} \label{lemThconfalg}
	Let $(M,g)$ be a 1-connected homogeneous plane wave. Then the Lie algebra $\conf(M,g)$ of conformal  vector fields of $(M,g)$ is a 1-dimensional extension of the Lie algebra of  Killing vector fields:
	$$\conf(M,g)=\Real D\oplus\isom(M,g)$$ and it consists of homothetic vector fields. 
\end{theorem}

\section{General lemmas}

In this section we prove 3 lemmas that will be used in the proof of Theorem~\ref{ThMain}.  

	\begin{lem}\label{lemexg}
	Let $(M=G/H,c)$ be a connected homogeneous pseudo-Riemannian conformal  manifold.
	Suppose that  a  Lie subgroup $\tilde{G}\subset G$  has  the open orbit  $U = \tilde{G}o  = \tilde{G}/\tilde{H}$,  where
	$  \tilde{H} \subset \tilde{G}$ is the  stability subgroup.
	If  the  isotropy group $j(\tilde{H})$   is a  subgroup of   the  orthogonal Lie group $\mathrm{O}(T_oU)$,  then  the  group $\tilde{G}$   preserves    the  metric  $g|_U$  which is  the  restriction  to  $U$ of  some   metric  $g \in c$ from the  conformal  class $c$.
	
\end{lem}

{\bf Proof.}   By the assumptions,   the value  $g_o$ at the point  $o$  of  any   metric $g \in c$ is invariant   with respect to the isotropy group  $j(\tilde{H})$.  Hence  it   can be extended   to  a $\tilde{G}$-invariant metric $g_U$ on  $U = \tilde{G}/\tilde{H}$. Since  $\tilde{G} $ is a subgroup   of the conformal group $G$,  the metric $g_U$    is  conformal to the   restriction  $g|_U$   to  $U$ of the  metric  $g \in c$. \qed

\begin{lem}\label{lemOpenOrbitTrans} Let $M  =  G/H$ be a connected  homogeneous manifold.
	Then any  normal   subgroup $F \subset G$  which   has  an open orbit   $U = Fo$  acts on  $M $ transitively.
	
\end{lem}

{\bf Proof.}  We remark  that    the   orbit  $Fx$   of any  point   $x \in  M$  is open. Indeed,  we may write
$x = ao$ for  some  $a \in  G$. Then $Fx = Fao = aFo$ is open.  This  implies  that   any orbit is  also closed.
Since  $M$ is  connected,   there is only one  orbit, i.e.  $F$  acts transitively on $M$. \qed

\begin{lem}\label{lempropF}    Let $(M=G/H, c)$ be a homogeneous  conformal pseudo-Riemannian manifold. 	Suppose that $F\subset G$ is a  normal Lie subgroup of $G$ acting transitively on $M$ by isometries of a metric $g\in c$. Then  $G$ acts by homothetic transformations of $g$.
\end{lem}

{\bf Proof.}
Let $a\in G$. There exists a function $\varphi$ such that $$a^*g=e^{2\varphi}g.$$ Let $f\in F$. Since the subgroup $F\subset G$ is normal, there exists an $f_1\in F$ such that $$af=f_1a.$$
Next, $$(af)^*g=f^*a^*g=f^*(e^{2\varphi}g)=e^{2f^*\varphi}g.$$
On the other hand,
$$(f_1a)^*g=a^*f_1^*g=a^*g=e^{2\varphi}g.$$
We conclude that $f^*\varphi=\varphi$. Since $F$ acts transitively on $M$, this implies that the function $\varphi$ is constant. \qed

A similar statement was proved in \cite{Bar1,Bar2}: if the identity connected component $F^0$ of the homothety group $F$ of a connected Lorentzian homogeneous complete manifold is normal in the conformal group $G$, then $F=G$.

\section{Proof of the Main Theorem}\label{secpMTh}

Let $(M=G/H, c)$ be a 1-connected non-conformally flat   essential homogeneous Lorentzian     conformal manifold   with   connected Lie   group  $G$.
Since $M$ is 1-connected and $G$ is connected, $H$ is connected.  
 Denote by $\g$ and $\h$ the Lie algebras of the Lie groups $G$ and $H$, respectively. Since $(M,c)$ is not conformally flat, by Corollary~\ref{Corconfflat}, the isotropy representation $$j:\h\to\co(V)=\Real \id_V \oplus \so(V)$$ in the tangent space 
 $$V=T_oM=\g/\h$$ is faithful. Hence we may identify the Lie algebra $\h$ with its image $j(\h)\subset\co(V)$, i.e.,
 $$\h\subset\co(V).$$
  Since the conformal group $G$ is essential and the stability subgroup $H$ is connected, we see that
$$\h\not\subset\so(V).$$

We  denote  by $$\tilde{\h} =  \h \cap \so(V)$$    the  codimension-one   ideal  
of  $\h$. Choose a  complementary    element  $$D = \id_V   +  C,\quad  C \in  \so(V)$$  such that
$$    \h  =  \Real  D  \oplus  \tilde{\h}.$$ The element $D$ is defined up to an element from $\tilde\h$. 


\subsection*{Step 1.} As the first step we  specify  the  element  $D$. 
	We prove the following 
\begin{prop}
	 Let $(M=G/H, c)$ be a   1-connected non-conformally flat  essential   homogeneous    conformal Lorentzian manifold. Then  the Lie  algebra  $\g$   admits  a   decomposition
	$$        \g  = \h \oplus V,\quad    V  =  \Real p  \oplus  E \oplus \Real q $$
	with the  stability    subalgebra $\h$   (identified  with the  isotropy  Lie  algebra $\h   =  j(\h) = \ad_{\h}|_V$)   of the  form
$$   \h =   \Real D  \oplus  \tilde{\h} ,\quad  \tilde{\h}\subset \so(V),\quad  D=\id_V-p\wedge q+C_0,\quad C_0\in\so(E).$$	\end{prop}

	\begin{lem}\label{lem1}
		Under the current assumptions it holds that $\id_V\not\in\h$.
			\end{lem}
	
	{\bf Proof.} Suppose that $\id_V\in\h$.
Let $\tilde \h=\h\cap\so(V)$. Then,
$$\h=\Real\id_V\oplus\tilde\h.$$
Since $[\id_{V}, \h]=0$ and $\ad_{\id_V}$ acts in $\g/\h$ as the identity,  there exists an $\ad_{\id_V}$-invariant subspace of $\g$ complementary to $\h$. This subspace may be identified with $V$. The decomposition $$\g=\h\oplus V$$ is the direct sum of eigenspaces of $\ad_{\id_V}$ corresponding to the eigenvalues $0$ and $1$.  From this and the Jacobi identity it follows that 
 $$[\tilde \h,V]\subset V,\quad  [V,V]=0.$$
Consider the subalgebra $$\tilde \g=\tilde \h \oplus V\subset \g.$$
Let $\tilde G\subset G$ be the connected Lie subgroup corresponding to the subalgebra $\tilde \g\subset\g$. The Lie subgroup $\tilde G \subset G$ is normal. Since $\tilde \g$ contains $V$, the $\tilde G$-orbit of the point $o$ is open. By Lemma \ref{lemOpenOrbitTrans}, $\tilde G$ acts transitively on $M$.
By Lemma \ref{lemexg}, there exists a metric $g\in c$ such that $\tilde G$ acts by isometries of $g$. The condition $[V,V]=0$ together with the equalities \eqref{NomizuMap} and \eqref{RoEq} given below imply that the metric $g$ is flat. This gives a contradiction. \qed


 Now  we  recall     a  description   of      endomorphisms  from the Lorentz
Lie  algebra  $\so(V)$.
It is  known    see, e.g., \cite{Geometry2}, that  there are   three types of   elements  $C$ of the  Lorentz  Lie  algebra  $\so(V)$,  described as  follows.
\begin{itemize}
	\item[{\bf Elliptic.}] 
	Such  element  $C =C_0$ annihilates  a  time-like   vector   $e_-\in V$ and  belongs  to the   orthogonal Lie  algebra
	$\so(E^{n+1}) \subset  \so(V)$ of the maximal Euclidean  subspace   $E^{n+1} = e_-^{\perp}$.
		\item[{\bf Hyperbolic.}] 
	With respect   to   some  Witt  basis   $p,e_1,\dots, e_n,q$ of $V$ it has  the  form 
	$$C=\alpha p\wedge q+C_0,\quad \alpha\in\Real,\quad \alpha< 0,\quad C_0\in\so(E),\quad E=\spa\{e_1,\dots,e_n\}.$$
	\item[{\bf Parabolic.}] 
		With respect   to   some  Witt  basis   $p,e_1,\dots, e_n,q$ of $V$ it may be written as $$C=\alpha p\wedge e_1+C_0,\quad \alpha\in\Real,\quad \alpha\neq 0,\quad C_0\in\so(E^{n-1}),\quad E^{n-1}=\mathrm{span}\{e_2,\dots e_n\}.$$
\end{itemize}

Let $C=\alpha p\wedge q+C_0$ be a hyperbolic element.  If we exchange the vectors $p$ and $q$, then $\alpha$ changes the sign, by that reason we assume that $\alpha<0$.

\begin{lem}\label{lem2} The Lie algebra $\h$ does not contain any element  $D=\id_V+C$ such that the element $C\in\so(V)$ is parabolic.
\end{lem}

{\bf Proof.} Suppose that $\h$ contains an element $D=\id_V+C$ such that 
$$C=\alpha p\wedge e_1+C_0\in\so(V)$$ is a parabolic element. Recall that $C_0\in\so(E)$ annihilates the vector $e_1$. The equality
$$D=(\id_V+C_0)+ \alpha p\wedge e_1$$ gives the decomposition of the element $D\in\co(V)$ into mutually commuting semisimple and nilpotent endomorphisms of $V$. Consequently the eigenvalues of $D$ acting on $V=\g/\h$ coincide with the eigenvalues of  $\id_V+C_0\in\co(V)$ and belong to the line $1+\Real i$. Similarly, from the decomposition \eqref{decompsoVC} it follows that the eigenvalues of $\ad_{D}$ on $\h\subset\co(V)$ belong to the set $\Real i$. Thus the eigenvalues of the endomorphism $\ad_D$ acting on $\g$ belong to the set $\Real i\cup (1+\Real i)$.  Consider the real  Jordan normal form \cite[Th. 3.4.1.2]{Horn}   of the endomorphism $\ad_D$ acting on $\g$. It is clear that the direct sum of  $\ad_{D}$-invariant subspaces of $\g$ corresponding to the eigenvalues from the set $\Real i$ coincides with $\h$, and the direct sum of  $\ad_{D}$-invariant subspaces of $\g$ corresponding to the eigenvalues from the set $1+\Real i$ is a vector subspace complementary to $\h$. 
 This vector subspace may be identified with the tangent space $V=\g/\h$. We obtain the   $\ad_{D}$-invariant decomposition $$\g=\h\oplus V.$$ 
The Jacobi identity
\begin{equation}\label{JacobiadD}\ad_{D}[X,Y]=[\ad_DX,Y]+[X,\ad_DY],\quad X,Y\in\g.\end{equation}
implies that if $\mathcal{A}$ and $\mathcal{B}$ are invariant subspaces of $\ad_D$ corresponding to the real Jordan blocks with the eigenvalues $\alpha$ and $\beta$, respectively, then $[\mathcal{A},\mathcal{B}]$ is contained in the $\ad_D$-invariant subspace corresponding to the eigenvalue $\alpha+\beta$. In particular, if $\alpha+\beta$ is not an eigenvalue, then 
$[\mathcal{A},\mathcal{B}]=0$.
This implies that $$[\h,V]\subset V,\quad [V,V]=0.$$ Thus $V\subset \g$ is a commutative ideal. As in the previous lemma, this implies that $(M,c)$ is conformally flat and we obtain a contradiction. \qed

	\begin{lem}\label{lem3} The Lie algebra $\h$ does not contain any 
		element $D=\id_V+C\in\h$ such that $C\in\so(V)$
		 is either  hyperbolic with $\alpha\neq - 1$ or elliptic.
	
\end{lem}

	{\bf Proof.} Suppose that $D=\id_V+C\in\h$, where $C$ is a  hyperbolic element defined by a number  $\alpha<0$ or elliptic; in the last case we assume that $\alpha=0$.
	It holds
	$$[D,p]=(1-\alpha)p,\quad [D,q]=(1+\alpha)q,\quad [D,E]=[\id_E+C_0,E]\subset E,$$
and	the eigenvalues of $D$ acting on $E$ belong to the set $1+ \Real i$. Hence the eigenvalues of $D$ acting on $V$ belong to the set
$$\{1\pm\alpha\}\cup (1+ \Real i).$$
 The eigenvalues of $D$ acting on $\co(V)$ belong to the set $$(\pm \alpha+\Real i)\cup\Real i.$$
We consider the real Jordan form of $\ad_D$ acting on $\g$, chose an $\ad_D$-invariant subspace of $\g$ complementary to $\h$ and identify it with the $D$-module $V$.

The eigenvalues of $D$ acting on $\wedge^2 V$ belong to the set
$$\{2\}\cup (2\pm\alpha+\Real i)\cup (2+ \Real i).$$
The map $$[\cdot,\cdot]|_{\wedge^2V}:\wedge^2V\to\h\subset\co(V)$$ is $D$-equivariant. 
This implies that if $$\alpha\not\in\left\{- 1,-2\right\},$$
then  $[V,V]=0$. As in the proof of Lemma \ref{lem2}, the equality $[V,V]=0$ implies that $(M,c)$ is conformally flat, which is not the case.

To prove the lemma it  remains to study the  case when   $C$   is a hyperbolic    element       with  $\alpha=-2$.  

{\bf Case    $\alpha = -2$.} 
The eigenvalues of $D$ acting on $V$ and $\h$ belong respectively to the sets
$$\{-1,3\}\cup (1+\Real i)\quad \text{and}\quad (\pm 2+\Real i)\cup \Real i.$$
The eigenvalues of $D$ acting on $\wedge^2 V$ and $\h\otimes V$ belong respectively to the sets
$$\{2\}\cup (4+\Real i)\cup (2+\Real i)\cup \Real i\quad\text{and}\quad (5+\Real i)\cup(\pm 3+\Real i)\cup(\pm 1+\Real i).$$
This implies that $$[\h,V]\subset V,\quad [V,V]\subset\h.$$
This means that $$\g=\h\oplus V$$ is a symmetric decomposition. Hence $(M,c)$ admits a locally symmetric Weyl connection
with the holonomy algebra $[V,V]\subset\co(V)$.
Results of \cite{DGS} show that any locally symmetric Weyl connection is closed, i.e., its holonomy algebra is contained in $\so(V)$. This means that $$[V,V]\subset\tilde\h.$$
The farther eigenvalues analysis of $D$ shows that $$[p,q]\in p\wedge E,\quad [E,q]\subset\Real p\wedge q\oplus \so(E),\quad [E,E]\subset p\wedge E.$$
If, for some $X\in E$, the projection of $[X,q]$ to
 $\Real p\wedge q$, is non-zero, then there exists an element $D'=\id+C_1\in\h$, where $C_1\in\so(E)$, which is impossible. Hence we may assume that
$$[E,q]\subset\so(E).$$
This and the Jacobi identity imply that $[p,q]=0$.
We see that there exist linear maps $$P:E\to\so(E),\quad Q:\wedge^2 E\to E$$  such that
$$[X,q]=P(X)\in \so(E),\quad [X,Y]=p\wedge Q(X,Y)\in p\wedge E,\quad X,Y\in E.$$ Let $E_1=Q(E,E)\subset E$. Denote by $E_2$  the orthogonal complement to $E_1$ in $E$.
The Jacobi identity implies
\begin{align}
	P(X)Y&=-Q(X,Y),\quad X\in E,Y\in E_1, \label{PQ1}\\
	Q(X,Y)+P(Y)X-P(X)Y&=0,\quad X,Y,Z\in E, \label{PQ2}\\
	(P(X)Y,Z)+(P(Y)Z,X)+(P(Z)X,Y)&=0,\quad X,Y,Z\in E. \label{PQ3}
	\end{align}	
From the last two equalities it follows that
$$(Q(X,Y),Z)=-(P(Z)X,Y),\quad X,Y,Z\in E.$$
Since $Q(X,Y)\in E_1$, the last equality implies that $P(E_2)=0$.
Let $X,Y\in E_1$. From \eqref{PQ1} and \eqref{PQ2} it follows that $P(X)Y=0$. Let $X\in E_1$ and $Y,Z\in E_2$.  From \eqref{PQ3} it follows that $P(X)Y=0$. Thus, $P=0$ and $Q=0$, i.e., $[V,V]=0$, and $(M,c)$ is conformally flat, which is a contradiction. 
 \qed

Lemmas    \ref{lem2} and \ref{lem3} imply   the  existence of      an  element $D\in\h$ of the form 
\begin{equation}\label{eqD}
D=\id_V-p\wedge q+C_0,\quad C_0\in\so(E).\end{equation}

\subsection*{Step 2.}

In the second step we will prove the following

\begin{prop}
	Let $(M=G/H, c)$ be a   1-connected non-conformally flat  essential   homogeneous    conformal Lorentzian manifold.
Then there exists an open neighbourhood  $U\subset M$ of the point $o$ and a metric $g\in c$ such that $(U,g|_U)$ is a  plane wave with the transitive action of  the isometry group of the metric $g|_U$.
\end{prop}

We have proved that there exists  $D\in \h$  given by \eqref{eqD}.
We claim that the element~$C_0\in \so(E)$ may be chosen in such a way that \begin{equation}\label{CondC0}
\Real C_0\cap (\tilde \h\cap \so(E))=\{0\},\quad\textrm{and} \quad [C_0,\tilde \h\cap \so(E)]\subset \tilde \h\cap \so(E).
\end{equation}
Indeed, for a chosen element $C_0$, let 
$$C_0=C_1+C_2,\quad C_1\in \tilde \h\cap \so(E),\quad C_2\in(\tilde \h\cap \so(E))^\bot,$$
where $(\tilde \h\cap \so(E))^\bot$ is the orthogonal complement to $\tilde \h\cap \so(E)$ in $\so(E)$. 
Since $$[C_0,\tilde \h\cap \so(E)]\subset \tilde \h\cap \so(E),\quad\textrm{and}\quad [\tilde \h\cap \so(E), (\tilde \h\cap \so(E))^\bot]\subset (\tilde \h\cap \so(E))^\bot,$$ the element $C_2$ satisfies conditions \eqref{CondC0}.
After the replacement of $D$ by $D-C_1$, $C_0$ will be replaced by $C_2$.
In what follows we assume that $C_0$ satisfies \eqref{CondC0}. This implies that
$$[D,\tilde \h\cap \so(E)]=0.$$ 

The eigenvalues of $\ad_D$ acting on $\Real p\wedge q\oplus\so(E)$, $p\wedge E$, $q\wedge E$ belong respectively to the sets $\Real i$, $1+\Real i$, $-1+\Real i$. Since $\ad_D$ preserves $\tilde\h$, this implies that
\begin{equation}\label{decomptildeh}
\tilde \h=\Big(\tilde \h\cap \big(\Real p\wedge q\oplus\so(E)\big)\Big)\oplus(\tilde \h\cap p\wedge E)\oplus(\tilde \h\cap q\wedge E).\end{equation}
If the projection of $\tilde \h\cap \big(\Real p\wedge q\oplus\so(E)\big)$ to $\Real p\wedge q$ is non-trivial, then we may change the element $D\in\h$ to an element $\id_V+C_1\in\h$, where $C_1\in\so(V)$ is elliptic, but this is impossible according to Lemma \ref{lem3}. Thus,
 \begin{equation}\label{decomptildeh1}
 \tilde \h=\big(\tilde \h\cap \so(E)\big)\oplus(\tilde \h\cap p\wedge E)\oplus(\tilde \h\cap q\wedge E)\end{equation}
is an $\ad_D$-invariant decomposition, and the eigenvalues of $\ad_D$ acting in the summands  belong respectively to the sets
$$\{0\},\quad 1+\Real i,\quad  -1+\Real i.$$
Next, $D$ preserves the decomposition 
$$V=\Real p\oplus\Real E\oplus\Real q,$$ and the corresponding eigenvalues belong respectively to the sets
$$\{2\},\quad 1+\Real i,\quad \{0\}.$$ 
 Let $\tilde p\in \g$ be an eigenvector of $\ad_D$ corresponding to the eigenvalue~$2$.  Let $$\g_0\subset\g\quad \text{and}\quad \mathcal{E}\subset \g$$ be the subspaces corresponding to real Jordan blocks with the eigenvalue $0$ and the eigenvalues from the sets  $1+\Real i$, respectively. Then
$$\Real D\oplus(\tilde\h\cap \so(E))\subset\g_0\quad\textrm{and}\quad \tilde \h\cap p\wedge E\subset \mathcal{E}.$$  Let 
$$\tilde E\subset  \mathcal{E}$$ be a complementary $\tilde\h\cap\so(E)$-invariant subspace. Let $\Real\tilde q\subset\g_0$ be an $\tilde\h\cap\so(E)$-invariant subspace
complementary to $\Real D\oplus(\tilde\h\cap \so(E))$.  We assume that the projection 
$$\g\to V=\g/\h$$ maps the vectors $\tilde p$ and $\tilde q$  respectively to the vectors $p$ and $q$.
We identify the subspace 
$\Real \tilde p\oplus\tilde E\oplus\Real \tilde q\subset\g$ with $V$.

The eigenvalues of $D$ acting on $\wedge^2 (\Real p\oplus \mathcal{E})$ belong to the set $$(2+\Real i)\cup (3+\Real i).$$ 
This shows that \begin{equation}\label{Eq*}
[\Real p\oplus \mathcal{E},\Real p\oplus \mathcal{E}]\subset  \Real p.\end{equation}
Next, it holds
$$[q,p]\subset\Real p,\quad [q,\mathcal{E}]\subset \mathcal{E}.$$
Let $E_1\subset E$ be the subspace such that $$p\wedge E_1=\tilde\h\cap p\wedge E.$$
We see that
$$\f=\mathcal{E}\oplus\Real p\oplus \Real q= p\wedge E_1\oplus V\subset \g$$
is a subalgebra. Since this subalgebra contains $V$, the orbit of the point $o$ for the corresponding connected Lie subgroup $F\subset G$ is an  open set $U$.  Denote by $F_o\subset F$ the stability subgroup of the point $o\in U$ under the action of $F$ on $U$. The corresponding Lie algebra is $\f_o=p\wedge E_1$.  By Lemma \ref{lemexg}, there exists a metric $g_U$ on $U$ belonging to $c_U$ such that $F$ is a transitive group of isometries of~$(U,g_U)$. 

\begin{lem}\label{lemhomplanewave} The homogeneous Lorentzian manifold $(U=F/F_o,g_U)$ is a homogeneous plane wave.
	\end{lem}

{\bf Proof.} As we have just seen, the Lie bracket of $\f$ satisfies
$$[q,p\wedge Y]=p\wedge KY-Y,\quad Y\in E_1,$$
where $K:E_1\to E_1$ is a linear map. If $K\neq 0$, then we change $V\subset\f$:
$$V\mapsto \Real p\oplus E'_1\oplus E_2\oplus\Real q,\quad E'_1=\{-p\wedge KY+Y|Y\in E_1\}.$$ Here  $E_2$ is the
orthogonal complement  to $E_1$ in $E$.
This allows us to assume that 
$$[q,p\wedge Y]=-Y,\quad Y\in E_1.$$
We conclude that 
$$\f=p\wedge E_1\oplus V=\f_o\oplus V$$ is a reductive decomposition of the Lie algebra $\f$.
The Lie bracket restricted to $V$ satisfies
\begin{align*}
[X,Z]&=\omega(X,Z)p,\\
[q,p]&=\lambda p,\\
[q,X]&=p\wedge BX+LX,
\end{align*}
where $X,Z\in E$, $Y\in E_1$, $\lambda\in\Real$, $B:E\to E_1$ and $L:E\to E$ are linear maps, and $\omega$ is a skew-symmetric form on $E$.

The Levi-Civita connection $\nabla$ of the metric   $g_U$    of the  reductive   Lorentzian  homogeneous  space    $(U=F/F_o,  g_U)$  is
determined by the Nomizu  operator   $$ \Lambda_V : V \to \so(V)$$  given by
\begin{equation}\label{NomizuMap}
  2 (\Lambda_V(X) Y,Z)  = ([X,Y]_V ,Z)  -  (Y, [X,Z]_V)-   (X, [Y,Z]_V), \quad X,Y,Z \in V, \end{equation}
where  $(.,.)$ is the  induced     Lorentz metric  in $V$ and      $X_V$  is  the  projection of a vector  $X \in  \g$ to  $V$, see, e.g., \cite{KNII}.
  If  $X^*$  is the   velocity vector  field of   1-parameter   group  $\exp (t X)$  generated  by an  element   $X \in V$,   then   \begin{equation}\label{nablaX*}\Lambda_V(X)  =  -\nabla X^*|_o.\end{equation} The  Nomizu operator      coincides (up to the sign) with the  covariant  derivative of an invariant   tensor field. The curvature tensor   at the point $o$ of the reductive homogeneous space 
$(U=F/F_o,  g_U)$  is given by
\begin{equation}\label{RoEq}R_o(X,Y)=[\Lambda_V(X),\Lambda_V(Y)]-\Lambda_V([X,Y]_V)-\ad([X,Y]_{\f_o}),\quad X,Y\in V,\end{equation} see, e.g., \cite{KNII}.
The Ambrose-Singer Theorem on holonomy for homogeneous spaces takes the form
$$\hol_o=\h_o+[\Lambda_V(V),\h_o]+[\Lambda_V(V),[\Lambda_V(V),\h_o]]+\cdots,$$
where $\h_o\subset\so(V)$ is the vector subspace spanned by the endomorphisms $R_o(X,Y)\in\so(V)$, $X,Y\in V$, see, e.g., \cite{KNII}.

It is easy to check that
\begin{align*}
\Lambda_V(p)&=0,\\
\Lambda_V(X)&=-\frac{1}{2}p\wedge\left(\omega +L+L^*\right)X,\quad X\in E,\\
\Lambda_V(q)&=\lambda p\wedge q+ \frac{1}{2}\left(-\omega +L-L^*\right),
\end{align*} where $\omega$ denotes the endomorphism of $E$ corresponding to the bilinear form $\omega$, and $L^*$ is the dual endomorphism to $L$ with respect to the scalar product on $E$. 
This immediately implies that the curvature tensor $R_o$ takes values in $p\wedge E\subset\so(V)$, and the holonomy algebra of the connection is contained in $p\wedge E\subset\so(V)$. Moreover,
\begin{equation}\label{eqlast1} R_o(X,Y)=0,\quad \forall X,Y\in p^\bot=\Real p+E.\end{equation}

Since the stability subgroup $F_o$ preserves the vector $p\in V$, the induced vector field $p^*$ on $U$ is $F$-invariant. The tensor field $\nabla p^*$ is  $F$-invariant as well. The equality $\Lambda_V(p)=0$ and \eqref{nablaX*} imply that 
\begin{equation}\label{eqlast4} \nabla p^*=0.\end{equation}
From \eqref{eqlast1} it follows that \begin{equation}\label{eqlast2} R(X,Y)=0,\quad \forall X,Y\in (p^*)^\bot.\end{equation}

The stability subgroup $\hat F_o$ of the $\hat F$-action on $U$ preserves the line $\Real p$ and the subspace $p^\bot\subset V$. Consequently, the distribution $(p^*)^\bot$ is  $\hat F$-invariant.
Recall that $\hat F$ preserves the Levi-Civita connection of the metric $g_U$ and the curvature tensor $R$. Hence the tensor field $\nabla_{|(p^*)^\bot} R$ is $\hat F$-invariant. This implies that the element $D$ annihilates the tensor $(\nabla_{|p^\bot} R)_o$. The tensor $(\nabla_{|p^\bot} R)_o$ may be considered as a linear map   
$$(\nabla_{|p^\bot} R)_o:(\Real p+E)\otimes \wedge^2 V\to p\wedge E.$$  The eigenvalues of $D$ acting on $p\wedge E$ belong to the set $1+\Real i$, and the eigenvalues of $D$-action on $(\Real p+E)\otimes \wedge^2 V$ do not belong to this set. This shows that $(\nabla_{|p^\bot} R)_o=0$. Hence, 
\begin{equation}\label{eqlast3} \nabla_{|(p^*)^\bot} R=0.\end{equation}
The equalities \eqref{eqlast4}, \eqref{eqlast2}, and \eqref{eqlast3}
show that $(U,g_U)$ is a plane wave. \qed

\subsection*{Step 3.} Now we are able to complete the proof of Theorem \ref{ThMain}.

Considering the elements from $\g$ as vector fields on $M$ and restricting them to $U$, we obtain the inclusion
$$\g \hookrightarrow\conf(U,g_U).$$ According to Theorem \ref{ThIsotrplw}, it holds
$$\conf(U,g_U)=\Real (\id_V-p\wedge q)\oplus\frakk\oplus p\wedge E\oplus V.$$
Consequently for the subalgebra $\tilde \h$ of the isotropy algebra $\h$
we obtain
$$\tilde\h\subset \frakk\oplus p\wedge E.$$
We conclude that
$$\tilde\h=\big(\tilde \h\cap \so(E)\big)\oplus (\tilde \h\cap p\wedge E).$$
This implies that the above defined subalgebra
$$\f=p\wedge E_1\oplus V=(\tilde\h\cap p\wedge E)\oplus V\subset \g$$ is an ideal, and the subgroup $F\subset G$ is normal. Consider, as at Step 2, the open orbit $U=Fo$ and the metric $g_U$ on it.  By Lemma \ref{lemOpenOrbitTrans}, $U=M$ and $g=g_U$ is a metric on $M$ from the conformal class $c$. By Lemma \ref{lemhomplanewave}, $(M,g)$ is a homogeneous plane wave. If $(M,g)$ is of type (a), then it is complete. If $(M,g)$ is of type (b), then, as we have seen in Section \ref{secHpw}, there exists a metric $g_1$ conformal to $g$ such that $(M,g_1)$ is a homogeneous plane wave of type (a), i.e., it is complete. Theorem \ref{ThMain} is now proved. \qed

\section{Conformal group of a homogeneous plane wave}\label{secCGHpvGroup}

In this section we  prove Theorem \ref{ThConf}.
Let $(M,g)$ be a 1-connected homogeneous  plane wave. Assume that $(M,g)$ is not conformally flat. First we suppose that $(M,g)$ is complete. Recall that a conformal transformation of a pseudo-Riemannian manifold is called a Liouville transformation if it preserves the Ricci tensor~\cite{KRLeovTrans}.

\begin{prop}
	  Let  $(M,g)$  be a  non  conformally    flat  complete  homogeneous plane  wave of  dimension  $n \geq 4$.
Then  each  conformal transformation  $a$ of $(M,g)$   is a Liouville transformation.
\end{prop}

{\bf Proof.} 
Fix an arbitrary point $x\in M$. 
Since the  isometry group of $(M,g)$ is transitive,   we may assume that a conformal transformation  $a$   preserves the point  $x \in M$.

Since $a$ is a conformal transformation, it preserves the Weyl curvature tensor. The curvature tensor $R$ of $(M,g)$ satisfies
$$R(\partial_{x^i},\partial_u)=\partial_v\wedge T\partial_{x^i},$$
where $$T=e^{uF}Be^{-uF},$$
and $R$ is zero on other pairs of basis vector fields. The Ricci tensor is given by
$$\Ric= \tr T (du)^2.$$
The structure of the Weyl curvature tensor is very similar to the structure of $R$,
 $$W(\partial_{x^i},\partial_u)=\partial_v\wedge \left(T-\frac{1}{n}{\textrm tr} T\,\id\right)\partial_{x^i}.$$
Consider the following Witt basis in $V=T_xM$:
$$p=(\partial_{v})_x,\quad e_i=(\partial_{x^i})_x,\quad q=(\partial_{v})_x-\frac{1}{2}g_x(\partial_u,\partial_u)(\partial_{u})_x.$$
The differential $a_x:T_xM\to T_x M$ belongs to the group
$$\CO(V)=\Real^* \id_V\cdot {\textrm O}(V).$$
There is a decomposition
$${\textrm O}(V)=\exp(q\wedge E)\cdot \Real^*\cdot {\textrm O}(E)\cdot \exp(p\wedge E),$$
where $$\Real^*\cdot {\textrm O}(E)=\left\{\left(
\begin{matrix}
\alpha&0&0\\
0&A&0\\
0&0&\alpha^{-1}
\end{matrix}\right),\alpha\in\Real^*,\,A\in {\textrm O}(E)\right\}.$$
The condition $a_xW_x=W_x$ implies
$$a_x\in \left\{\alpha\id_V\cdot\left(
\begin{matrix}
\alpha&0&0\\
0&A&0\\
0&0&\alpha^{-1}
\end{matrix}\right),\alpha\in\Real^*,\,A\in H\right\}\cdot\exp(p\wedge E),$$ where $H\subset {\textrm O}(E)$ is the subgroup commuting with $T$. This immediately implies that $a_x$ preserves~$\Ric_x$. \qed


Since the metric $g$ is complete, according to \cite[Corollory 1]{KRLeovTrans},  each Liouville transformation of $(M,g)$ is a homothety transformation. This completes the proof of Theorem \ref{ThConf} for homogeneous plane waves of type (a). 

As we heave seen in Section \ref{secHpw}, each homogeneous plane wave of type (b) is of the form $(M,e^{u}g)$, where $(M,g)$ is a homogeneous plane wave of type (a),  and $u$ is the global coordinate as in Section \ref{secHpw}. The groups of conformal transformations of  $(M,e^{u}g)$ and $(M,g)$ coincide. It is clear that the group of conformal transformations of $(M,g)$ is generated by the group of isometries and by the transformations \eqref{homoth}. According to \cite[Section 6]{HMZ}, each isometry of $(M,g)$ transforms the coordinate $u$ to $u+c$ for a constant $c$. This shows that each isometry of $(M,g)$ is a homothety transformation of $(M,e^{u}g)$. Likewise, each transformation \eqref{homoth} of $(M,g)$ is a homothety transformation of  $(M,e^{u}g)$. This completes the proof of Theorem \ref{ThConf}. 
\qed

\section{Special cases}\label{secEx}

\subsection*{Lie groups with conformal vector fields induced by derivations}

 In \cite{ZC2021,ZC2024}, the authors studied simply connected  Lie groups $Q$ with
pseudo-Riemannian left-invariant metrics $g$ that admit  essential conformal transformations induced by  derivations of the corresponding Lie algebras. A classification is obtained in the Lorentzian signature.

Let us apply our results to this situation. Let $(Q,g)$ be a simply connected Lie group with a left-invariant Lorentzian metric $g$.
Denote by $V$ the Lie algebra of $Q$. Suppose that $D$ is a derivative of $V$ that induces an essential conformal transformation of $(Q,g)$. Consider the Lie algebra
$$\g=\Real D\oplus V,$$ where the Lie bracket of $D$ and $V$ is given by the action of $D$ on $V$. By the assumption,
$$\g\subset\conf(Q,g).$$
Let $G$ be the corresponding connected Lie subgroup of the Lie group of conformal transformations of $(Q,g)$. The group $G$ acts transitively as an essential group of conformal transformations of
$(Q,g)$. It is clear that the isotropy representation is faithful. By Theorem \ref{ThMain}, $(Q,g)$ is a homogeneous plane wave, and $G$ acts on  $(Q,g)$ by homothetic transformations. There exists a decomposition
$$V=\Real p\oplus E\oplus \Real q$$ such that $D$ is given by \eqref{eqD}. The eigenvalues analysis of $D$ shows that the Lie bracket of $V$ satisfies
\begin{align*}
[X,Z]&=\omega(X,Z)p,\quad [p,X]=0,\\
[q,p]&=\lambda p,\\
[q,X]&=LX,
\end{align*}
where $X,Z\in E$ are arbitrary,  $\lambda\in\Real$, and $L:E\to E$ is linear maps.
The Jacobi identity is equivalent to the following conditions:
\begin{equation}\label{omegaL}\omega(LX,Y)+\omega(X,LY)=\lambda\omega(X,Y),\quad X,Y\in E,\end{equation}
$$[C_0,L]=0,\quad C_0\omega=0.$$
The inclusion $$\g\hookrightarrow\conf(Q,g)$$ defines the element $D\in\conf(Q,g)_o$ and the $D$-invariant subspace $$V\subset \conf(Q,g)$$ complementary to $\conf(Q,g)_o$. We obtain a decomposition
$$\conf(Q,g)=\Real D\oplus\frakk\oplus p\wedge E\oplus V.$$
It holds $$[q, p\wedge X]=-X+p\wedge KX,\quad X\in E,$$
where $K$ is an endomorphism of $E$.

The eigenvalues of $D$ acting on $p\wedge E$ and $E$ are the same, so instead of  $E\subset \g$ we may consider any subspace of the form
\begin{equation}\label{E'}E'=\{X+p\wedge \varphi X|X\in E\},\end{equation} where $\varphi$ is an endomorphism of $E$ commuting with $D$ and $\tilde \h\cap \so(E)$.
It holds \begin{align}\label{newomega}[X+p\wedge \varphi X,Y+p\wedge \varphi Y]&=\omega(X,Y)p+(X,\varphi Y)p-(Y,\varphi X)p,\\
\label{newL}[q,X+p\wedge \varphi X]&=(L-\varphi)X+p\wedge \varphi(L-\varphi)X+p\wedge  (K\varphi-\varphi(L-\varphi))X.\end{align}
For an endomorphism $A$ of $E$ we denote by $A^s$ and $A^{sk}$ the symmetric and skew-symmetric parts of $A$.
Let \begin{equation}\label{phi}
\varphi=\frac{1}{2}\omega+L^s,\end{equation} where $\omega$ is the endomorphism of $E$ identified with the bilinear form $\omega$.
We see that $E$ may be chosen in such a way that
\begin{equation}\label{condE}\omega=0,\quad L \text { is skew-symmetric}.\end{equation}
Moreover, for a fixed vector $q$,  the condition \eqref{condE} determines the subspace $E\subset \g$ uniquely.
Let us denote $L$ by $F$.
   If $\lambda\neq 0$, then we apply the change
$$p\mapsto \lambda p,\quad q\mapsto \frac{1}{\lambda}q.$$
This allows us to assume that $\lambda$ is either 0 or 1.
Now the Lie bracket on $p\wedge E\oplus V\subset\conf(Q,g)$ is  exactly as the Lie bracket of $(M,g)$ from Section \ref{secHpw} with
\begin{align}
\label{F} F&=L^{sk}-\frac{1}{2}\omega,\\
\label{B} B&=\lambda L^s+[L^{sk},L^s]-\frac{1}{4}\omega^2-(L^s)^2.
\end{align}
Conversely, it holds

\begin{prop} Let $(M,g)$ be a 1-connected homogeneous plane wave given by $\lambda,F,B$. Then   $(M,g)$ is a Lie group with a left-invariant metric $g$ if and only if there exist $\omega$ and $L$ satisfying \eqref{omegaL}, \eqref{F}, and \eqref{B}.\end{prop}

\subsection*{Cahen-Wallach spaces}
A 1-connected Cahen-Wallach space $(M,g)$ may be characterized as homogeneous plane wave of type (a) with $F=0$. Thus $(M,g)$ is determined by the symmetric endomorphism $B$. A Cahen-Wallach space $(M,g)$ is a symmetric space, and  isometry Lie algebra 
$\isom(M,g)$ admits the $\mathbb{Z}_2$-grading 
\begin{equation}\label{algisomCW}\isom(M,g)=(\frakk\oplus p\wedge E)\oplus V,\quad
V=\Real p\oplus E\oplus \Real q, 
\end{equation} 
where $\frakk\subset\so(V)$ is the subalgebra commuting with $B$. 
The Lie bracket of  $\isom(M,g)$ restricted to $V$ coincides up to the sign with the curvature tensor of $(M,g)$, and it holds 
$$[q,p]=[p,X]=[X,Y]=0,\quad
[q,X]=p\wedge BX, \quad X,Y\in E.$$ 

 The conformal transformations of Cahen-Wallach spaces have been studied in details in \cite{LT2022}. In particular, it is explained that if
$(M,g)$ is not conformally flat, then the conformal transformations of $(M,g)$ are homothetic transformations, and the
group of conformal transformations of $(M,g)$ is a semidirect product of $\Real$ and the isometry group of $(M,g)$.
Theorem~\ref{ThConf} is a generalization of this result.

Here we use our results to study the question, when there is a Lie group acting simply transitively on a Cahen-Wallach space.
From the above we see that a simply connected Cahen-Wallach space $(M,g)$ is a Lie group with a left-invariant Lorentzian metric if and only if there exists an endomorphism $L$ such that
\begin{equation}\label{CWLg}L^{sk}L^s+L^sL^{sk}=0, \quad B=[L^{sk},L^s]-(L^{sk})^2-(L^s)^2.\end{equation}
Decomposing $E$ as the direct sum of $\ker L^s$ and $(\ker L^s)^\bot$, we see that the first condition in \eqref{CWLg} is equivalent to the condition $L^{sk}|_{(\ker L^s)^\bot}=0$, or, in other words, the images of $L^{sk}$ and $L^{s}$ are orthogonal.
Consequently,
\begin{equation}\label{BLskLs}B=-(L^{sk})^2-(L^s)^2.\end{equation} It is clear that the eigenvalues of $-(L^s)^2$ are non-negative, while the eigenvalues of $-(L^{sk})^2$ are positive and the multiplicity of each eigenvalue of $-(L^{sk})^2$ is even.

We have proved

\begin{prop} Let $(M,g)$ be a 1-connected Cahen-Wallach space defined by a symmetric endomorphism~$B$ of $E$. Then   $(M,g)$ is a Lie group with a left-invariant metric $g$ if and only if each positive eigenvalue of $B$ has even multiplicity.
	\end{prop}

\begin{prop} Let $(M,g)$ be a 1-connected Cahen-Wallach space defined by a symmetric endomorphism~$B$ of $E$. 
	Suppose that $(M,g)$ admits a simply transitive action of a Lie group $Q$ of isometries, i.e., there exists an endomorphism $L$ of $E$ such that \eqref{BLskLs} holds and the images of $L^{sk}$ and $L^{s}$ are orthogonal. Then the Lie algebra $V$ of the Lie group $G$ is isomorphic to
	$$V'=\Real p\oplus E\oplus \Real q$$ with the Lie bracket  
$$[X,Z]=2(L^{sk}X,Z)p,\quad [p,X]=[q,p]=0,\quad [q,X]=LX.$$
The inclusion $$V'\hookrightarrow \isom(M,g)$$ is given by
$$p\mapsto p,\quad X\mapsto X+p\wedge L X,\quad q\mapsto q.$$
\end{prop}

Let $(Q,g)$ be a simply connected Lie group with a left-invariant Lorentzian metric $g$ and suppose as above that $(Q,g)$ admits an essential conformal transformation induced by a derivation of the corresponding Lie algebra. The metric $g$ is bi-invariant if and only if the scalar product on $V$ is $\ad_V$-invariant. It is easy to see that this condition is equivalent to the equalities $\lambda=0$, $\omega=L^{sk}=0$. This implies the following known result.
\begin{prop} \cite{KO2019} A 1-connected Cahen-Wallach space $(M,g)$ is a Lie group with a bi-invariant Lorentzian metric if and only if the eigenvalues of $B$ are non-positive.\end{prop}


\begin{thebibliography}{10}

\bibitem{A2017} 	
D. Alekseevsky,
Lorentzian manifolds with transitive conformal group.
 Note di Matematica 37 (2017), 35--47.

\bibitem{Geometry2} D. V. Alekseevsky, E. B. Vinberg, and A. S. Solodovnikov, Geometry of spaces of constant curvature. Geometry II, Encyclopaedia Math. Sci., vol. 29, Springer, Berlin, 1993,  1--138.

\bibitem{A}  D.V. Alekseevsky, Groups of  conformal  transformations  of Riemannian
spaces.  Math. Sb.  89 (1972), no. 1, 280--296.

\bibitem{A1} D.V.  Alekseevsky, Self-similar Lorentzian manifolds.  Ann. Global
Anal. Geom. 3 (1985),  no. 1, 59--74.

\bibitem{A2} D.V. Alekseevsky, The sphere and the Euclidean space are the only Riemannian manifolds with essential conformal transformations. Uspekhi Math. Nauk 28 (1973),  no. 5, 289--290.




\bibitem{Blau} M. Blau, Plane Waves and Penrose Limits. Lecture Notes http://www.blau.itp.unibe.ch/lecturesPP.pdf

\bibitem{BO} M. Blau, M. O'Loughlin. Homogeneous plane waves. Nuclear Physics B 654 (2003), 135--176.

\bibitem{Bar1} C. Barbance, 
Transformations conformes des variétés lorentziennes homogenes. 
C. R. Acad. Sci., Paris, Sér. A 291 (1980), 342-350.


\bibitem{Bar2} C. Barbance, Transformations conformes des variétés lorentziennes homogenes.
Tensor, New Ser. 39 (1982), 173-178.

\bibitem{HMZ}
M. Hanounah, L. Mehidi, A. Zeghib, On homogeneous plane waves. J. Math. Phys. 66(5) (2025).

\bibitem{HSII2024}
J. Holland, G. Sparling, Sachs equations and plane waves II: Isometries and conformal isometries. arXiv:2405.12748

\bibitem{Horn} R.A. Horn,  C.R. Johnson,   Matrix Analysis (2nd ed.). Cambridge Univ. Press, 2021.

\bibitem{DGS} A. Dikarev, A.S. Galaev, E. Schneider, Recurrent Lorentzian Weyl spaces. J. Geom. Anal.   34, 282 (2024).



\bibitem{Fer}  J. Ferrand, The action of conformal transformations on a
Riemannian manifold.  Math. Ann. 304 (1996), no. 2, 277--291.

\bibitem{F}  C. Frances, Sur les vari\'et\'es lorentziennes dont le groupe
conforme est essentiel,
Math. Ann. 332 (2005), no. 1, 103--119.

\bibitem{F1}  C. Frances, Essential conformal structures in
Riemannian  and  Lorentzian structures, in "Recent Development of
Pseudo-Riemannian geometry",  ed. D.V.Alekseevsky, H. Baum,
ESI Lect. Math. Phys., Eur. Math. Soc., Z\"urich, 2008, 234--260.



\bibitem{F-Z} C. Frances, A. Zeghib, Some remarks on pseudo-Riemannian
conformal actions of simple Lie groups. Math. Res. Lett. 12 (2005), 49--56.

\bibitem{KO2019} I. Kath and M. Olbrich. Compact quotients of Cahen-Wallach spaces. Mem. Amer. Math. Soc., 262(1264), v+84,
2019.

\bibitem{KNII} S. Kobayashi, K. Nomizu, Foundations of Differential Geometry, Vol. 2, Wiley and Sons, New York, 1969.

\bibitem{KRLeovTrans} W. K\"uhnel and H.-B. Rademacher, Conformal diffeomorphisms preserving the
Ricci tensor. Proc. Amer. Math. Soc. 123 (1995), 2841--2848.




\bibitem{LT2022} T. Leistner,  S. Teisseire,  Conformal transformations of Cahen-Wallach spaces. Annales de l'Institut Fourier (2025).



\bibitem{Obata} M. Obata, The conjectures on conformal
transformations of Riemannian manifolds. J . Diff. Geom.
6 (1971), 247--258.



\bibitem{P} M.N. Podoksenov,  Conformally  homogeneous  Lorentzian
manifolds, Sib. Mat. J. 33 (1992), no. 6,   154--161.

\bibitem{ZC2021}
H. Zhang, Z. Chen,
Lie groups with conformal vector fields induced by derivations.			Journal of Algebra
584 (2021), 304-316.


\bibitem{ZC2024}
H. Zhang, Z. Chen,
On Lie groups with conformal vector fields induced by derivations.		Transformation Groups (2024). 



\end{thebibliography}
\end{document}